\documentclass[12pt]{amsart}
\usepackage{amssymb}

\newtheorem{thm}{Theorem}[section]
\newtheorem{prop}[thm]{Proposition}
\newtheorem{lemma}[thm]{Lemma}

\newtheorem{remark}[thm]{Remark}

\newtheorem{example}[thm]{Example}

\numberwithin{equation}{section}

\newcommand{\Edot}{{E_\bullet}}

\newcommand{\Gr}{{\bf G}}

\newcommand{\G}{{\mathbb{G}}}
\newcommand{\K}{{\mathbb{K}}}

\begin{document}

\title[Transversality in the Schubert calculus]{Elementary transversality
        in the Schubert calculus  in any Characteristic}   

\author{Frank Sottile}
\address{Department of Mathematics and Statistics\\
        University of Massachusetts\\
        Amherst, Massachusetts 01003\\
        USA}
\email{sottile@math.umass.edu}
\urladdr{http://www.math.umass.edu/\~{}sottile}
\date{6 July 2002}
\subjclass{14N10, 14M15, 14N15, 14N35, 14P99}
\thanks{Research supported in part by NSF grant DMS-0070494}
\keywords{enumerative geometry, transversality, Schubert variety} 

\begin{abstract}
We give a characteristic-free proof
that general codimen\-sion-1 Schubert varieties meet transversally in a
Grassmannian and in some related varieties. 
Thus the corresponding intersection numbers
are enumerative in all characteristics. 
Existing transversality results do not apply to these enumerative
problems, emphasizing the need for additional theoretical work on
transversality.
We also strengthen some results in enumerative real algebraic 
geometry.
\end{abstract} 

\maketitle


\section*{Introduction}

Schubert~\cite{Sch1879} declared enumerative geometry to be concerned with
all questions of the following form:
How many geometric figures of some type satisfy certain given conditions?
For conditions imposed by general fixed figures, 
the approach---then, as now---is to interpret the conditions as
subvarieties of a parameter space of figures, which give cycle classes
in the Chow ring of that parameter space.
Then the degree of the product of these cycle classes provides an algebraic
count of the solutions, counting the solutions weighted with certain
intersection multiplicities. 
Thus this degree solves the original problem of enumeration---if each
solution occurs with multiplicity one, that is, if the subvarieties meet
transversally at each point of intersection.

This demonstrates how the validity of the standard approach to enumerative
geometry via multiplying cycles classes in a Chow ring rest upon the following
basic premise in enumerative geometry:
General subvarieties of the parameter space meet transversally at each point
of intersection.
An intersection number is enumerative if this basic premise holds for the
corresponding intersection.
Kleiman~\cite{MR50:13063} established this basic premise in characteristic
zero, when an algebraic groups acts transitively on the parameter space.
The restriction to characteristic zero is necessary.
General translates of arbitrary
subvarieties simply do not meet transversally in positive characteristic.
Kleiman~\cite{MR50:13063} exhibits a subvariety of a grassmannian
that does not meet any translate of a particular codimen\-sion-2 Schubert variety
transversally. 

However, Kleiman's example does not arise in an enumerative geometric problem.
In fact, in every known case, general Schubert subvarieties of flag varieties
meet generically transversally (transverse along a dense subset of each
component), in any characteristic.
Thus the question remains:
To what extent does this basic principle of enumerative geometry hold in
positive characteristic?
Laksov and Speiser develop a general theory for transversality---they
give a condition, using tangent spaces to families of subvarieties, that
implies a general member of a family meets any fixed
subvariety transversally~\cite{La75,LaSp90,LaSp92,Sp88}. 
By Kleiman's example, families of codimen\-sion-2 Schubert subvarieties do not
satisfy this condition.

By Theorem~E of~\cite{So97a}, general Schubert subvarieties of a grassmannian
of 2-planes in a vector space meet transversally, in any characteristic. 
Here, we give an elementary and characteristic-free 
proof that general simple (codimen\-sion-1 or divisorial) Schubert varieties meet 
transversally, when the ambient space is one of the following: 
(i) the grassmannian, (ii) the flag manifold for the special linear group, 
(iii) the orthogonal grassmannian, or
(iv) the space of parameterized rational curves of a fixed degree in a
grassmannian.
While these results are rather special in that they only involve simple
Schubert varieties, such enumerative problems are in fact quite natural.
More importantly, these results suggest that transversality is
ubiquitous in enumerative geometry.

Thus, the corresponding intersection numbers are enumerative for fields of any
characteristic (except characteristic 2 for the orthogonal grassmannian).
This includes some genus zero Gromov-Witten invariants of the grassmannian, by
(iv).
For example, given 12 general 4-planes in 7-dimensional space,
there are exactly 462 3-planes that meet all 12.
(This number was computed by Schubert~\cite{Sch1886c}.)
Similarly, given $N=5q+6$ general points $s_1,s_2,\ldots,s_N$ in 
${\mathbb P}^1$ and $N$ general 3-planes $K_1,K_2,\ldots,K_N$ in 5-dimensional
space, there are exactly $F_{5+5q}$ (the $(5+5q)$th Fibonacci number) 
degree $q$ maps  $M$ from ${\mathbb P}^1$ to the grassmannian of 2-planes in
5-space satisfying $M(s_i)\cap K_i\neq\{0\}$, for each
$i=1,\dotsc,N$~\cite{In91,RRW96}.

We  use the formalism of~\cite{So00c} to show transversality.
While that was developed to construct real-number solutions to
these enumerative problems, it can also be used to show transversality.
We state our results in Section 1, review the formalism of~\cite{So00c}
in Section 2,  and prove our results in Section 3.
In Section 4, we strengthen the real enumerative
geometric results of~\cite{So99a,So00,So00c}.
Finally, in Section 5, we show that families of
simple Schubert varieties do not satisfy the condition of 
Laksov and Speiser.

While the simple Schubert varieties we study are members
of a linear series of sections of an ample line bundle, the general section
is not a Schubert variety, and so standard Bertini-type theorems in positive
characteristic do not apply.
Our results, together with the inapplicability of the general theory of Laksov 
and Speiser, emphasize the need for  new, perhaps less general or
more specific, criteria that imply transversality.
For example, is there a Bertini-type theorem for the Grassmannian concerning
intersections with a general codimen\-sion 1 Schubert variety?

Similarly, our results and extensive computer calculations
(for example, see~\cite[\S4.3]{So_M2}) suggest that the classical Schubert
calculus of enumerative geometry (and the quantum Schubert calculus) is
enumerative in all characteristics (except characteristic 2 for the orthogonal
groups).
Specifically, we make the following conjecture, for which we feel there is
ample support (both theoretical and examples).
\medskip

\noindent{\bf Conjecture. }
{\it 
Let $\K$ be an algebraically closed field and $X$ a flag variety for
a reductive algebraic group defined over $\K$.
If $X_1,X_2,\ldots,X_s$ are Schubert varieties of $X$ in general position,
the sum of whose codimensions equals the dimension of $X$, then \vspace{-4pt}
the intersection
$$
  X_1\cap X_2\cap\cdots\cap X_s\vspace{-2pt}
$$
is transverse.
}\medskip

We thank Dan Laksov for his thoughtful and constructive comments on an earlier
version of this manuscript.


\section{Statement of results}
Let $\K$ be any infinite field.
We describe the spaces and simple Schubert varieties for which we prove
transversality. 

\begin{enumerate}
\item[(i)] 
 Let $0<r<n$ be integers.
 Let ${\bf G}(r,n)$ denote the grassmannian of $r$-planes in 
 $\K^n$.
 An $(n{-}r)$-dimensional subspace ($(n{-}r)$-plane) $K$ defines a simple
 (co\-di\-men\-sion-1) Schubert subvariety of ${\bf G}(r,n)$
$$
  \Omega(K)\ :=\ \{H\in{\bf G}(r,n)\mid H\cap K\neq\{0\}\}\,.
$$

\item[(ii)]
 For a sequence ${\bf r}:0<r_1<r_2<\cdots<r_m<n$ of integers, let 
 ${\mathbb F}\ell_{\bf r}$ be the variety of $m$-step flags
$$
  \Edot\ :\ E_{1}\subset E_{2}\subset\cdots\subset
            E_{m}\subset\K^n
$$
 with $\dim E_{i}=r_i$ for $i=1,\ldots, m$.
 An $(n{-}r_i)$-plane $K$ defines a simple Schubert subvariety of 
 ${\mathbb F}\ell_{\bf r}$
$$
  \Phi_i(K)\ :=\ \{\Edot\in{\mathbb F}\ell_{\bf r}\mid
                   E_{r_i}\cap K\neq\{0\}\}\,.
$$

\item[(iii)]
 Suppose that the characteristic of $\K$ is not 2 and
 $n=2r{+}1$ is odd.
 Let $\langle\,\ ,\ \rangle$ be a non-degenerate split symmetric
 bilinear form on $\K^n$, that is, one for which there is a basis 
 $e_1,\ldots,e_n$ of $\K^n$ such that 
 \begin{equation}\label{form}
  \langle e_i,\,e_j\rangle\ =\ \left\{
    \begin{array}{lcl} 1&\ &\mbox{if }i+j=n+1=2r+2\\
                       0&  &\mbox{otherwise}
    \end{array}\right..
 \end{equation}
 The resulting quadratic form is
 \begin{equation}\label{quadratic}
  q\Bigl(\sum x_ie_i\Bigr)\ =\ 2x_1x_{2n-1}\,+\,2x_2 x_{2n-2}\,+\,\cdots
                   \,+\,2x_rx_{r+2}\,+\,x_{r+1}^2\,.
 \end{equation}
 A subspace $H$ of $\K^n$ is {\it isotropic} if 
 $\langle\cdot,\cdot\rangle|_H\equiv 0$, and isotropic subspaces have
 dimension at most $r$.
 The orthogonal grassmannian $OG(r)$ is the space of maximal (dimension-$r$)
 isotropic subspaces of $(\K^n,\langle\cdot,\cdot\rangle)$.
 This has dimension $\binom{r+1}{2}$ and, by our choice of the
 form~\eqref{form}, the $\K$-rational points are dense.
 (This may be deduced, for example from the coordinates for Schubert cells
 of the orthogonal flag manifold given in~\cite[p.~67]{FP98}.)
 An isotropic $r$-plane $K$ of $\K^n$ defines a simple Schubert
 subvariety of $OG(r)$
$$
  \Psi(K)\ :=\ \{H\in OG(r)\mid H\cap K\neq\{0\}\}\,.
$$

\item[(iv)]
 For an integer $q\geq 0$, let ${\mathcal M}^q_{r,n}$ be the space of degree
 $q$ maps $M\colon{\mathbb P}^1\to{\bf G}(r,n)$.
 A point $s\in{\mathbb P}^1$ and an $(n{-}r)$-plane $K$ define a simple
 (quantum) Schubert subvariety of ${\mathcal M}^q_{r,n}$
 $$
  Z(s,K)\ :=\ \{M\in{\mathcal M}^q_{r,n}\mid
                      M(s)\cap K\neq \{0\}\}\,.
 $$
\end{enumerate}

Each of these spaces have more general \emph{Schubert varieties}.
We will prove the following theorem.

\begin{thm}\label{thm:main}
 Suppose $\K$ is infinite.
 Let $X$ be either ${\bf G}(r,n)$, ${\mathbb F}\ell_{\bf r}$,
 $OG(r)$, or ${\mathcal M}^q_{r,n}$ and suppose $Z\subset X$ is any Schubert
 variety of $X$.
 Then general simple Schubert varieties $Y_1,Y_2,\ldots,Y_{\dim Z}$ meet
 $Z$ transversally on $X$.
\end{thm}

We prove Theorem~\ref{thm:main} by exhibiting 1-parameter families 
of simple Schubert varieties having properties which
imply that general members of these families intersect
transversally. 
We remark that in each case except the classical flag manifold, there is only
one simple Schubert variety, up to the action of the relevant group
($GL_n(\K)$, $O_{2n+1}(\K)$, or $GL_n(\K)\times PGL_2(\K)$).
However, the theory we use involves the action of the 1-dimensional torus
$\G_m$ and the families we use are parameterized by $\G_m$, so our arguments
really do involve different families. 

\section{The Method of Schubert Induction}

We review the method of Schubert induction introduced in~\cite{So00c}.
Let $\G_m$ be the group scheme whose $\K$-rational points
are the invertible elements $\K^\times$ of $\K$,
considered as a group under multiplication.
While this theory holds for families of subvarieties over any curve, we use
$\G_m$ as the base for our families.
Suppose that $X$ is a projective variety and we have a subvariety 
$E\subset X\times\G_m$, where the fibres $E\to\G_m$ are equidimensional, that
is, a family $E\to\G_m$ of subvarieties of $X$.
Then the scheme-theoretic limit $\lim_{s\to 0}E(s)$ is defined to be the fibre
over 0 of the closure $\overline{E}$ of $E$ in $X\times\K$.
Since $\G_m$ is a curve, this fibre has the expected dimension.
(See Remark 9.8.1 of~\cite{Ha77}.)

A {\it Bruhat decomposition} of an irreducible variety $X$
defined over a field $\K$ 
is a finite decomposition
$$
  X\ =\ \coprod_{w\in I} X^\circ_w
$$
satisfying the following conditions.
\begin{enumerate}
\item[(1)]
     Each stratum $X^\circ_w$ is a locally closed irreducible
     subvariety defined over $\K$.
\item[(2)]
     The closure $\overline{X^\circ_w}$ of a stratum is a union of some
     strata $X^\circ_v$.  
\item[(3)]
     There is a unique 0-dimensional stratum $X^\circ_{\hat{0}}$.
\end{enumerate}
For $w\in I$, define the {\it Schubert variety} $X_w$ to be
$\overline{X^\circ_w}$. 
The space $X$ is the top-dimensional Schubert variety.
By condition (2), the intersection $X_w\cap X_v$ of two Schubert
varieties is a union of some Schubert varieties.
The {\it Bruhat order} on $I$ is the order induced by inclusion of Schubert
varieties: $u\leq v$ if $X_u \subseteq X_v$.
Set $|w|:=\dim X_w$.
 
Let ${\mathcal Y}\rightarrow \G_m$ be a family of
divisors on $X$. 
For $s\in\G_m$, let $Y(s)$ be the fibre of ${\mathcal Y}$ over the
point $s$.
We say that ${\mathcal Y}$ {\it respects} the Bruhat decomposition if, for
every  $w\in I$, the (scheme-theoretic) limit 
$Z:=\lim_{s\rightarrow 0} ( Y(s)\cap X_w)$ is supported on a union of Schubert
subvarieties $X_v$ of codimension $1$ in $X_w$.
Write $v\prec\!\!_{\mathcal Y}\,\,w$ when $X_v$ is the support of a
component of the limit scheme $Z$.
This relation generates a suborder of the Bruhat order, and, as
we explain below, the number of solutions to an enumerative
problem given by simple Schubert varieties equals the number of saturated
chains in such a suborder.
Our notation for this suborder, $\prec\!\!_{\mathcal Y}$, reflects its
dependence on the family $\mathcal{Y}$.
A family $\mathcal{Y}\to\G_m$ of divisors of $X$ is {\it multiplicity-free}
if it respects the Bruhat 
decomposition and if each component of the scheme-theoretic limit $Z$ is
reduced at its generic point. 

The general fibre $Y(s)$ of such a multiplicity-free family of divisors on $X$
meets each Schubert variety $X_w$ generical\-ly trans\-ver\-sal\-ly.
Indeed, if the general fibre of $\mathcal{Y}$ did not meet $X_w$ generically
transversally, then every intersection $Y(s)\cap X_w$ would have a non-reduced
component.
This would imply that $\lim_{s\to 0}Y(s)\cap X_w$ has a non-reduced component,
violating our assumption that the family $\mathcal{Y}$ is multiplicity-free.

Thus when $X$ is smooth and $Y(s)$ is a general member of the family
$\mathcal{Y}$, we have the formula   
\begin{equation}\label{eq:cycle_prod}
  [X_w]\cdot[Y(s)]\ =\ \  [X_w \cap Y(s)]\ =\ 
    \sum_{v\prec_{{\mathcal Y}}w} \, [X_v]\,
\end{equation}
in the Chow ring of $X$. 
The second equality in~\eqref{eq:cycle_prod} expresses the rational
equivalence of the fibres of the family $\mathcal{Y}\cap(X_w\times\K)$ and the
first equality is a basic property of any intersection theory, namely that a
generically transverse intersection of subvarieties of $X$ 
represents the product of their cycle classes, as the intersection
multiplicities are equal to 1~\cite[Remark 8.2, p.~138]{Fu98}. 

A collection ${\mathcal Y}_1,\ldots,{\mathcal Y}_l$ of families of
divisors of $X$ meets the Bruhat decomposition of $X$
properly\footnote{The term `in general position with respect to the Bruhat
decomposition' is used in~\cite{So00c}.} 
if for all $w\in I$, for general $s_1,s_2,\ldots,s_l\in\G_m$,  
for each $1\leq k\leq l$, and for every $k$-subset $\{i_1,i_2,\ldots,i_k\}$ of  
$\{1,2,\ldots,l\}$, the intersection
 \begin{equation}\label{eq:intersection}
   Y_{i_1}(s_{i_1}) \cap Y_{i_2}(s_{i_2})\cap \cdots \cap 
   Y_{i_k}(s_{i_k})\cap X_w
 \end{equation}
is {\it proper} in that either it is empty or else it has (the expected)
dimension  $|w| - k$.

Suppose that we have a collection 
${\mathcal Y}_1,{\mathcal Y}_2,\ldots,{\mathcal Y}_l$  
($l=\dim X$) of multiplicity-free families of divisors of $X$ with each
meeting the Bruhat decomposition of $X$ properly.
Recall that $\hat{0}$ is the index of the minimal Schubert variety, which is a
simgle point.
For $w\in I$, let $\deg(w)$ count the number of (saturated) chains in the
Bruhat order
 \begin{equation}\label{eq:chains}
    \hat{0}\prec_1 w_1\prec_2 w_2\prec_3\cdots\prec_{k-1} w_{k-1}\prec_k
     w_k = w\,,
 \end{equation}
where $\prec_i = \prec_{{\mathcal Y}_i}$ and $|w|=k$.
This is the degree of  the intersection
 \begin{equation}\label{E:G-Int}
   Y_{1}(s_{1}) \cap Y_{2}(s_{2})\cap \cdots \cap 
   Y_{k}(s_{k})\cap X_w\,.
 \end{equation}
A result of~\cite{So00c} asserts that this intersection is
free of multiplicities.

\begin{prop}[\cite{So00c}]\label{prop:schubert_induction}
 Suppose $X$ has a Bruhat decomposition, 
 ${\mathcal Y}_1,{\mathcal Y}_2,\ldots,{\mathcal Y}_l$ are  
 multi\-pli\-ci\-ty-free families of divisors in $X$ over $\G_m$ meeting the
 Bruhat decomposition of $X$ properly, and each respects the Bruhat
 decomposition.
 Then, for every $1\leq k\leq l$ and every $w\in I$ with $|w|=k$,
 the intersection~\eqref{E:G-Int} is transverse for general
 $s_1,s_2,\ldots,s_k\in\G_m$ and has degree $\deg(w)$. 
 In particular, when $\K$ is algebraically closed, such an
 intersection consists of $\deg(w)$ reduced points.
\end{prop}

The point of Proposition~\ref{prop:schubert_induction} is that while we only
assume that the intersections~\eqref{eq:intersection} have the expected
dimension, we are able to conclude that~\eqref{E:G-Int} is transverse and to
compute the number of solutions without reference to the Chow ring.
There is a simple criterion (Lemma 2.4 of~\cite{So00c}) which implies that
a collection of families meets the Bruhat decomposition properly.
 
\begin{lemma}\label{lem:2.4}
 Suppose a variety $X$ has a Bruhat decomposition.
 Let ${\mathcal Y}_1,{\mathcal Y}_2,\ldots,{\mathcal Y}_l$ be a collection
 of families of divisors of $X$.
 If each family ${\mathcal Y}_k\to\G_m$ satisfies
$$
  \bigcap_{s\in\G_m}{\mathcal Y}_k(s)\ =\ \emptyset\,,
$$
 then the collection of families 
 ${\mathcal Y}_1,{\mathcal Y}_2,\ldots,{\mathcal Y}_l$ meets the Bruhat
 decomposition properly.
\end{lemma}

The elementary methods of~\cite{So00c} are illustrated by the 
following  proof.
If the families ${\mathcal Y}_1,{\mathcal Y}_2,\ldots,{\mathcal Y}_l$ do not
meet the Bruhat decomposition properly, then after possibly reindexing, there
is a $w\in I$, integer $k$, and 
points $s_1,s_2,\ldots,s_k\in\G_m$ such that
 \begin{equation}\label{eq:lemma}
   Y_{1}(s_{1}) \cap Y_{2}(s_{2})\cap \cdots \cap 
   Y_{k}(s_{k})\cap X_w
 \end{equation}
has (the expected) dimension of $|w|-k\geq 0$, but for every 
$s\in\G_m$, the intersection 
$$
   Y_{1}(s_{1}) \cap Y_{2}(s_{2})\cap \cdots \cap 
   Y_{k}(s_{k})\cap  Y_{k+1}(s)\cap X_w
$$
also has dimension $|w|-k$.
But then some component of~(\ref{eq:lemma}) lies in every
subvariety $Y_{k+1}(s)$, contradicting the assumption of the lemma.
\qed


\section{Proof of Theorem~1.1}
By Proposition~\ref{prop:schubert_induction}, for each space ${\bf G}(r,n)$,
${\mathbb F}\ell_{\bf r}$, $OG(r)$, and ${\mathcal M}^q_{r,n}$, we need only
to construct multiplicity-free families of simple  Schubert
subvarieties over $\G_m$ such that 
the entire collection meets the Bruhat decomposition properly.
For the space ${\mathcal M}^q_{r,n}$ of rational curves in a grassmannian,
we work in Drinfeld's compactification ${\mathcal K}^q_{r,n}$,
also called the quantum grassmannian.

\subsection{The grassmannian}
Let $e_1,e_2,\ldots,e_n$ be an ordered basis for $\K^n$.
For a sequence $\alpha:1\leq\alpha_1<\alpha_2<\cdots<\alpha_r\leq n$, the 
 Schubert variety $\Omega_\alpha$ is
$$
  \Omega_\alpha\ :=\ \{H\in\Gr(r,n)\mid \dim H\cap F_{\alpha_j}\geq j\ 
     \mbox{ for }j=1,2,\ldots, r\}\,,
$$
where $F_i$ is the linear span of $e_1,e_2,\ldots, e_i$.
Set $\Omega^\circ_\alpha:=
 \Omega_\alpha-\{\Omega_\beta\mid\beta\leq\alpha\}$.
Here, $\leq$ is given by componentwise comparison, $\alpha\leq\beta$ if and
only if $\alpha_i\leq\beta_i$ for all $i$.
Write $\binom{[n]}{r}$ for the resulting partial order on this set of
sequences. 
With these definitions, the grassmannian ${\bf G}(r,n)$ has a Bruhat decomposition 
$$
  \Gr(r,n)\ =\ \coprod_{\alpha\in\binom{[n]}{r}} \Omega^\circ_\alpha\,,
$$
indexed by sequences $\alpha\in\binom{[n]}{r}$.
For $\alpha\in\binom{[n]}{r}$, set $|\alpha|=\sum(\alpha_i-i)$, which is the
dimension of $\Omega_\alpha$.
Write $\beta\lessdot \alpha$ when $\beta<\alpha$ with $|\beta|=|\alpha|-1$.

If $H$ is the row space of an $r$ by $n$ matrix,
then the $\binom{n}{r}$ maximal minors of that matrix are the 
Pl\"ucker coordinates of $H$.
Write these as $p_\alpha$ for $\alpha\in\binom{[n]}{r}$.
These Pl\"ucker coordinates give an embedding of $\Gr(r,n)$ into 
${\mathbb P}^{\binom{n}{r}-1}$.

We use the following elementary set-theoretic fact,
originally due to Schubert~\cite{Sch1886c}.
 \begin{equation}\label{E:Simple-Pieri}
    \Omega_\alpha\bigcap\{p_\alpha{=}0\}\ =\ 
        \bigcup_{\beta\lessdot\alpha}\Omega_\beta\,,
 \end{equation}
and the intersection is generically transverse.
This is a consequence of the related ideal- (or scheme-) theoretic fact,
which may be deduced from the combinatorics of the Pl\"ucker ideal of 
${\bf G}(r,n)$.
There are, however, very elementary reasons that 
intersection~\eqref{E:Simple-Pieri} is as claimed (set-theoretically) and is
generically transverse.
Indeed, if we set 
\[
   K\ :=\ \mbox{Span}\langle e_i\mid i\in\{1,2,\dotsc,n\}-
            \{\alpha_1,\dotsc,\alpha_r\} \rangle\,,
\]
then $\Omega(K)$ has equation $p_\alpha=0$ in the Pl\"ucker coordinates, and a
simple computation in local coordinates proves the equality
in~\eqref{E:Simple-Pieri}, as well as showing that it is generically transverse.
(For a synthetic argument, see Theorem~2.4(2) of~\cite{So97d}.)\medskip

\noindent{\it Proof of Theorem~$1.1$ for ${\bf G}(r,n)$.}
 Let $K\subset\K^n$ be an $(n{-}r)$-plane, none of whose
 Pl\"ucker coordinates vanish.
 Since no Pl\"ucker coordinate vanishes identically on the grassmannian, such
 a plane exists as $\K$ is infinite.
 Let $\G_m$ act on $\K^n$ by  $s.(e_j) = s^je_j$ and set $K(s):=s.K$.
 Let ${\mathcal Y}(K)\to\G_m$ be the family of simple Schubert
 varieties whose fibre over $s\in\G_m$ is the Schubert variety $\Omega(K(s))$.

 Theorem \ref{thm:main} is a consequence of
 Proposition~\ref{prop:schubert_induction} and the 
 following claim.

\noindent{\bf Claim.}
   Let $K_1,K_2,\ldots,K_l$ be $(n{-}r)$-planes, each with no vanishing
  Pl\"ucker coordinates.
\begin{enumerate}
\item[(i)] 
      Each family ${\mathcal Y}(K_i)$ preserves the Bruhat decomposition of  
      ${\bf G}(r,n)$ and is multi\-pli\-city-free.
\item[(ii)] 
      The collection of families
      ${\mathcal Y}(K_1),{\mathcal Y}(K_2),\ldots,{\mathcal Y}(K_l)$ 
      meets the Bruhat decomposition properly.
\end{enumerate}

\noindent{\it Proof of Claim. }
 Represent $K_i$ as the row space of an $(n{-}r)$ by $n$ matrix $K_i$.
 Then $K_i(s)$ is represented by the same matrix, but with column $j$
 multiplied by $s^j$. 
 If a $r$-plane $H$ is the row space of an $r$ by $n$ matrix $H$, then 
$$
  H\cap K_i(s)\ \neq\ \{0\}\qquad \Longleftrightarrow\qquad 
   \det\left[ \begin{array}{c} K_i(s)\\ H\end{array}\right]\ =\ 0\,.
$$
 Laplace expansion of the determinant along the rows of $H$ gives the
 equation in Pl\"ucker coordinates  for $H$ to lie in
 $\Omega(K_i(s))$
\begin{equation}\label{eq:plucker_equation}
  0 \ =\ \sum_{\beta\in\binom{[n]}{r}} p_\beta k_\beta 
           s^{r(n-r)+\binom{n-r+1}{2}-{|\beta|}}\,,
\end{equation}
 where $k_\beta$ is the appropriately signed maximal minor of $K_i$ for the
 columns complementary to $\beta$.
 Up to a sign, this is a Pl\"ucker coordinate of $K_i$.

 If we restrict~(\ref{eq:plucker_equation}) to 
 the Schubert variety $\Omega_\alpha$ and divide by the common factor
 $s^{r(n-r)+\binom{n-r+1}{2}-{|\alpha|}}$, the result is a polynomial 
 with constant term $p_\alpha k_\alpha$.
 Thus the limit scheme 
 $\lim_{s\to 0}(\Omega_\alpha\cap \Omega(K_i(s)))$
 is defined in $\Omega_\alpha$ by $p_\alpha=0$, and so
 by~\eqref{E:Simple-Pieri}, 
$$ 
  \lim_{s\to 0}(\Omega_\alpha\cap \Omega(K_i(s)))\ =\ 
     \sum_{\beta\lessdot \alpha} \Omega_\beta\,,
$$
 which proves the first claim.

 For a fixed $H\in{\bf G}(r,n)$,~(\ref{eq:plucker_equation}) is a non-zero
 polynomial in $s$, and so there are finitely many 
 $s\in \G_m$ with $H\in\Omega(K_i(s))$.
 In particular, $\bigcap_{s\in\G_m}\Omega(K_i(s))=\emptyset$.
 By Lemma~\ref{lem:2.4}, this implies the second claim and completes the
 proof of Theorem~\ref{thm:main} for ${\bf G}(r,n)$.
\qed

\subsection{The flag manifold}
Fix an ordered basis
$e_1,e_2,\ldots,e_n$ for $\K^n$.
Let  ${\bf r}:0<r_1<r_2<\cdots<r_m<n$ be integers.
The flag manifold ${\mathbb F}\ell_{\bf r}$ has a Bruhat decomposition
$$
  {\mathbb F}\ell_{\bf r}\ =\ \coprod X^\circ_w\,,
$$
indexed by those permutations $w=w_1w_2\ldots w_n$ in the symmetric group
${\mathcal S}_n$ on $n$ letters whose descent set
$\{i\mid w_i>w_{i+1}\}$ is a subset of $\{r_1,r_2,\ldots,r_m\}$.
Here, $X^\circ_w$ is a Schubert cell of ${\mathbb F}\ell_{\bf r}$
(see \S 9 of~\cite{Fu97}).

Let $\G_m$ act on $\K^n$ by $s.e_i=s^ie_i$.
Suppose $K$ is an $(n{-}r_i)$-plane, none of whose Pl\"ucker coordinates
vanish.
For $s\in\G_m$, set $K(s):=s.K$.
Let ${\mathcal Y}(K)\to\G_m$ be the family whose fibre over 
$s\in\G_m$ is the simple Schubert variety
$\Phi_i(K(s))$.
The projection ${\mathbb F}\ell_{\bf r}\to{\bf G}(r_i,n)$  sends a 
flag $\Edot\in{\mathbb F}\ell_{\bf r}$ to its $i$th component
$E_i\in{\bf G}(r_i,n)$, and $\Phi_i(K(s))$ is the inverse image of
$\Omega(K(s))$.
It is an easy consequence (see~\cite[\S 2]{So00c} for details) of the
results proven in \S 3.1 that
the family ${\mathcal Y}(K)$ preserves the Bruhat decomposition of 
the flag variety ${\mathbb F}\ell_{\bf r}$ and is multiplicity-free.
(This last fact follows from Monk's formula~\cite{Monk}.)
Since $\bigcap_{s\in\G_m}\Omega(K(s))=\emptyset$, 
we have $\bigcap_{s\in\G_m}\Phi_i(K(s))=\emptyset$, 
and so we conclude that if we have 
any $i_1,i_2,\ldots,i_l\in\{1,2,\ldots,m\}$ and 
$(n-r_{i_j})$-planes $K_j$ for $j=1,\dotsc,l$, 
none of whose Pl\"ucker coordinates vanish, then the families 
${\mathcal Y}(K_1),{\mathcal Y}(K_2),\ldots,{\mathcal Y}(K_l)$
meet the Bruhat decomposition properly.

These facts together establish Theorem~\ref{thm:main} for the flag manifold
${\mathbb F}\ell_{\bf r}$.\qed

\subsection{The orthogonal grassmannian}
Suppose that the characteristic of $\K$ is not 2 and $n=2r+1$ is odd.
Let $e_1,e_2,\ldots,e_n$ be a basis for $\K^n$
for which the symmetric bilinear form is as given by~(\ref{form}).
The orthogonal grassmannian has a Bruhat decomposition
$$
  OG(r)\ =\ \coprod X^\circ_\lambda\,,
$$
indexed by  decreasing sequences $\lambda$ of integers
$n\geq\lambda_1>\lambda_2>\cdots>\lambda_l>0$.
This decomposition is induced from that of
${\bf G}(r,n)$ by the inclusion 
$\iota\colon OG(r)\hookrightarrow {\bf G}(r,n)$.
As before, $X^\circ_\lambda$ is a Schubert cell of $OG(r)$.

Let $\G_m$ act on $\K^n$ by
$s.e_i=s^ie_i$.
This induces an action of $\G_m$ on the
orthogonal grassmannian, as the quadratic form~\eqref{quadratic} is
homogeneous of degree $2n+2$ under this action.
Given $K\in OG(r)$ with no vanishing Pl\"ucker
coordinates, set $K(s):=s.K$ and
let ${\mathcal Y}(K)\to\G_m$ be the family
of simple Schubert varieties with fibres 
$\Psi(K(s))$ ($=\iota^{-1}(\Omega(K(s))$, set-theoretically).

As in \S 3.2 (see~\cite[\S 3]{So00c}), the results proven for the family
$\Omega(K(s))$ in  \S 3.1 imply the corresponding results for the family
${\mathcal Y}(K)$, and any collection of such
families.
(Multiplicity-freeness follows from a cohomological formula due to
Chevalley~\cite{CH91}.) 

In this way, we  establish Theorem~\ref{thm:main} for the orthogonal
grassmannian $OG(r)$.\qed

\subsection{The space of rational curves in the grassmannian}

The space ${\mathcal M}^q_{r,n}$ of degree $q$ maps
$M\colon{\mathbb P}^1\to{\bf G}(r,n)$ is a smooth quasi-projective 
variety of dimension $qn+r(n{-}r)$.
The Pl\"ucker coordinates of such a map are homogeneous forms of degree
$q$.
Choosing $\K\subset{\mathbb P}^1$,
these forms are polynomials of degree $q$ in the parameter
$s\in\K$.
Let $z_{\alpha^{(a)}}$ be the coefficient of $s^{q-a}$ in the $\alpha$th
Pl\"ucker coordinate of a map $M$.
These coefficients give quantum Pl\"ucker coordinates for 
${\mathcal M}^q_{r,n}$, determining the Pl\"ucker embedding of
${\mathcal M}^q_{r,n}$ into the projective space
${\mathbb P}(\bigwedge^p\K^{m+p}\otimes\K^{q+1})$.
The closure of ${\mathcal M}^q_{r,n}$ in this embedding
is the singular Drinfeld compactification or quantum grassmannian 
$\mathcal{K}^q_{r,n}$. 

Let ${\mathcal C}^q_{r,n}:=\{\alpha^{(a)}\mid\alpha\in\binom{[n]}{r}
    \ \mbox{and}\ 0\leq a\leq q\}$ be the set of indices of
quantum Pl\"ucker coordinates.
This set is partially ordered as follows
$$
  \alpha^{(a)}\ \leq\ \beta^{(b)} \quad\Longleftrightarrow\quad
   a\leq b \, \mbox{ and }  \,\,
   \alpha_i\leq \beta_{b-a+i} \,
  \mbox{ for }\, i = 1,2,\ldots,p-b+a. 
$$
The quantum Schubert varieties
$$
  \overline{Z_{\alpha^{(a)}}}\ :=\ \{z\in{\mathcal K}^q_{r,n}\mid 
    z_{\beta^{(b)}}=0\ \mbox{ if }\ \beta^{(b)}\not\leq \alpha^{(a)} \}\ ,
$$
are the Schubert varieties of a Bruhat decomposition of 
${\mathcal K}^q_{r,n}$~\cite{So00}
$$
  {\mathcal K}^q_{r,n}\ =\ \coprod_{\alpha^{(a)}\in{\mathcal C}^q_{r,n}}
          Z^\circ_{\alpha^{(a)}}\,.
$$
Here $Z^\circ_{\alpha^{(a)}}$ is the set of points in 
$\overline{Z_{\alpha^{(a)}}}$ with non-vanishing coordinate $z_{\alpha^{(a)}}$.

Let $I^q_{r,n}$ be the ideal of the quantum grassmannian.
The ideals of quantum Schubert varieties have a very simple description, 
which generalizes~\eqref{E:Simple-Pieri}.

\begin{prop}[\cite{SS_SAGBI}]\label{q-ideals}
\mbox{ }
\begin{enumerate}
\item[(i)] The ideal $I_{\alpha^{(a)}}$ of \,$\overline{Z_{\alpha^{(a)}}}$ is
           $I^q_{r,n}+\langle z_{\beta^{(b)}}\mid 
             \beta^{(b)}\not\leq\alpha^{(a)}\rangle$.
\item[(ii)] ${\displaystyle \langle 
          I_{\alpha^{(a)}},z_{\alpha^{(a)}}\rangle\ =\ 
         \bigcap_{\beta^{(b)}\lessdot \alpha^{(a)}} I_{\beta^{(b)}}}$.
\end{enumerate}
\end{prop}

This was established modulo embedded primes of lower dimension (which is
sufficient for our purposes) in~\cite{RRW96,RRW98}.\medskip

\noindent{\it Proof of Theorem~$1.1$ for ${\mathcal M}^q_{r,n}$.}
Let $K\subset\K^n$ be an $(n{-}r)$-plane, none of whose
Pl\"ucker coordinates vanish.
Let $\G_m$ act on $\K^n$ by 
$s.(e_i) = s^ie_i$ and set $K(s):=s.K$.
Consider the family of simple Schubert varieties of
${\mathcal M}^q_{r,n}$
 \begin{equation}\label{sqsv}
  Z(s,K)\ :=\ \{M\in{\mathcal M}^q_{r,n}\mid
                      M(s^n)\cap K(s)\neq \{0\}\}\,.
 \end{equation}

Let ${\mathcal Y}(K)\to\G_m$ be the family of subvarieties of
${\mathcal K}^q_{r,n}$ whose fibre $\overline{Z(s,K)}$ over $s\in\G_m$ is the 
closure of $Z(s,K)$.
As in~\cite[\S 3]{So00},  expanding the determinantal equation 
for $M$ to lie in $Z(s,K)$ gives
the linear equation for this fibre 
 \begin{equation}\label{qeq}
   0\ =\ \sum_{\alpha^{(a)}\in{\mathcal C}^q_{r,n}}
                z_{\alpha^{(a)}} k_\alpha 
     s^{qn+r(n-r) + \binom{n-r+1}{2}-|\alpha^{(a)}|}\,.
 \end{equation}

As in \S 3.1, the form of this equation and 
Proposition~\ref{q-ideals} show that the family ${\mathcal Y}(K)$ respects
the Bruhat decomposition of ${\mathcal K}^q_{r,n}$ and is multiplicity-free.
Furthermore, given any $(n{-}r)$-planes $K_1,K_2,\ldots,K_l$ in 
$\K^n$, none of whose Pl\"ucker coordinates vanish, the resulting
families ${\mathcal Y}(K_1)$,
${\mathcal Y}(K_2),\ldots,{\mathcal Y}(K_l)$ meet the Bruhat decomposition
properly. 
Thus general members of these families  meet transversally on
the quantum grassmannian ${\mathcal K}^q_{r,n}$, and hence on its dense
subset ${\mathcal M}^q_{r,n}$.
This completes the proof of Theorem~\ref{thm:main}.
\qed\medskip

In general, all points of intersection of general members
of the families ${\mathcal Y}(K_1)$,
${\mathcal Y}(K_2),\ldots,{\mathcal Y}(K_{\dim{\mathcal M}^q_{r,n}})$ lie in
the space  ${\mathcal M}^q_{r,n}$ of curves.
An argument given in~\cite{So00} uses work of Bertram~\cite{Be97} concerning
the quot scheme compactification of 
${\mathcal M}^q_{r,n}$.
This argument is valid here, as the pertinent results of Bertram hold for over 
arbitrary fields.
(See~\cite{So_Lowell} for a survey of this quantum intersection problem.)


\section{Some more reality} 
We strengthen the results of~\cite{So99a,So00,So00c}:

\begin{prop}\label{prop:real}
 Suppose $\K={\mathbb R}$.
 Let $X$ be one of the spaces ${\bf G}(r,n)$, ${\mathbb F}\ell_{\bf r}$,
 $OG(r)$, or ${\mathcal M}^q_{r,n}$.
 There exist simple Schubert varieties $Y_1,Y_2,\ldots,Y_{\dim X}$ that meet
 transversally in the complexification $X_{\mathbb C}$ of $X$, and every
 point of intersection is real.
\end{prop}

In~\cite{So99a,So00,So00c}, we constructed families of simple Schubert
varieties defined by subspaces $K(s)$ which 
osculate a given rational normal curve.
Unlike the families constructed in \S 3, those families respect the
Bruhat order only in characteristic zero.
The calculations of \S 3 enable a more flexible choice of subspaces.

An action of the torus ${\mathbb R}^\times$ on a real vector space $V$ of
dimension $n$ is {\it general} if $V$ is a direct sum of
1-dimensional eigenspaces, each with a different character.
Given such an action, we say that a linear subspace $L\subset V$ is 
{\it general} if it none of its Pl\"ucker coordinates vanishes, where we
define Pl\"ucker coordinates with respect to a basis of eigenvectors.
When $n=2r{+}1$ and $V$ is equipped with a non-degenerate (split)\footnote{We
require the form to be split so that $OG(r)$ has sufficiently many real
points.} 
symmetric bilinear form of signature $\pm1$, we require the
form to be homogeneous with respect to this action.
See~\eqref{form} for an example of such a split form for the general action
$s.e_i=s^i\cdot e_i$.
A general action induces actions of the torus on
the spaces ${\bf G}(r,n)$, ${\mathbb F}\ell_{\bf r}$, and $OG(r)$.
We obtain a torus action on ${\mathcal M}^q_{r,n}$ by 
having ${\mathbb R}^\times$ act on the source ${\mathbb P}^1$ of the
maps via $s.[a,b]:=[a,s^Nb]$, for some integer $N$.
This action on ${\mathbb P}^1$ will be general (and induce a general action
on ${\mathcal M}^q_{r,n}$) when $N>N_0$, for some $N_0$
(described below) depending upon the given general torus action.

\begin{thm}\label{thm:action}
 Let $V$ be an $n$-dimensional real vector space equipped with a general
 action of the torus and equip ${\mathbb P}^1$ with a general action
 of the torus.
 Let $X$ be one of the spaces ${\bf G}(r,n)$, ${\mathbb F}\ell_{\bf r}$,
 $OG(r)$, or ${\mathcal M}^q_{r,n}$.
 For any simple Schubert varieties $Y_1,Y_2,\ldots,Y_{\dim X}$ defined by
 general linear subspaces of $V$, there exist real numbers
 $s_1,s_2,\ldots,s_{\dim X}$ such that the translates
 $s_1.Y_1,s_2.Y_2,\ldots,s_{\dim X}.Y_{\dim X}$ meet transversally
 in the complexification $X_{\mathbb C}$ of $X$, and every
 point of intersection is real.
\end{thm}

There is a second part to Proposition~\ref{prop:schubert_induction}
of~\cite{So00c}, which we use for the proof below.

\begin{prop}[\cite{So00c}]\label{reality}
 Under the same hypotheses as Proposition~\ref{prop:schubert_induction},
 but with $\K={\mathbb R}$,
 there exist points $s_1,s_2,\ldots,s_l\in {\mathbb R}$ such that
 for every $1\leq k\leq l$ and every $w\in I$ with $|w|=k$,
 the intersection
$$
   Y_{1}(s_{1}) \cap Y_{2}(s_{2})\cap \cdots \cap 
   Y_{k}(s_{k})\cap X_w
$$
 is transverse and each of its $\deg(w)$ points are real.
\end{prop}

\noindent{\it Proof of Theorem~\ref{thm:action}. }
The characters of ${\mathbb R}^\times$ are the monomials $s^i$ for $i$
an integer.
We use the integer $i$ to represent the character $s^i$.
Suppose $V$ is an $n$-dimensional real vector space equipped with a general
action of ${\mathbb R}^\times$.
Assume that $e_1,e_2,\ldots,e_n$ is a basis of eigenvectors of $V$ with
respective characters $i_1,i_2,\ldots,i_n$ where
$i_1<i_2<\cdots<i_n$.
Then the arguments of \S\S 3.1, 3.2, and 3.3 remain valid with this
action in place of the action $s.e_i=s^ie_i$.
This action induces the action on Pl\"ucker coordinates
$s.p_\alpha=s^{J(\alpha)}p_\alpha$, where 
$J(\alpha):=\sum_j i_{\alpha_j}$.
The key facts are that the map $\alpha\mapsto J(\alpha)$ is an order preserving
map from the poset $\binom{[n]}{r}$ to the integers, and 
that the exponent 
$r(n{-}r)+\binom{n-r+1}{2}-{|\alpha|}$ of~(\ref{eq:plucker_equation})
is replaced by $i_{r+1}+\cdots+i_n-J(\alpha)$.
The theorem follows from Proposition~\ref{reality}, for the spaces ${\bf G}(r,n)$, 
${\mathbb F}\ell_{\bf r}$, and  $OG(r)$.

For the space of rational curves ${\mathcal M}^q_{r,n}$, 
set $N_0:=i_n{-}i_1{+}1$ and suppose $N>N_0$ for the torus action on 
${\mathbb P}^1$.
Thus the action of $s\in{\mathbb R}^\times$ on the map $M$
is given by the map $t\mapsto  s.[M(s^Nt)]$.
On the quantum Pl\"ucker coordinates, this is 
$s.z_{\alpha^{(a)}}= s^{I(\alpha^{(a)})}z_{\alpha^{(a)}}$,
where $I(\alpha^{(a)})= Na + J(\alpha)$.
Again the map $\alpha^{(a)}\mapsto I(\alpha^{(a)})$ is
an order preserving map from the poset 
${\mathcal C}^q_{r,n}$ to the integers.
The arguments of \S 3.4 extend to this more general action,
a key point being the exponent
$qn+r(n-r) + \binom{n-r+1}{2}-|\alpha^{(a)}|$  of~(\ref{qeq}) is now 
replaced by 
$qN+i_{r+1}+\cdots+i_n-I(\alpha^{(a)})$.
Thus the theorem follows for the space ${\mathcal M}^q_{r,n}$, by
Proposition~\ref{reality}.
\qed\medskip

The proof of Proposition~\ref{reality} in~\cite{So00c}
(see~\cite[\S 4]{So00}) gives further
information about the choice of 
the points $s_1,s_2,\ldots,s_l$.
By $\forall s_1\gg s_2\gg\cdots\gg s_l$, we mean 
$$
  \forall s_1>0\ \ \exists \epsilon_2>0,\ {\rm such\ that}\ 
  \forall s_2<\epsilon_2\ \  \cdots\ \ \exists \epsilon_l>0,\ 
  {\rm such\ that} \ \forall  s_l<\epsilon_l\, .
$$
Thus the existential statement ``there exist points
$s_1,s_2,\ldots,s_{\dim X}\in{\mathbb R}$'' of Theorem~\ref{thm:action}
may be replaced by 
``$\forall s_1\gg s_2\gg\cdots\gg s_{\dim X}$''.

\section{Not a determinantal pair}
Laksov and Speiser~\cite{La75,Sp88,LaSp90,LaSp92}
develop the notion of a determinantal family of
subvarieties thereby giving a general criterion for proving that a general
member of a family meets arbitrary subvarieties transversally.
This approach studies the tangent spaces of members of the family
with respect to possible tangent spaces for arbitrary subvarieties.
To facilitate applications of their theory, they give a local 
condition which implies that a family is determinantal.
We show that this local condition does not hold for families of simple
Schubert varieties in ${\bf G}(r,n)$, and thus does not
hold for the other families of simple Schubert
varieties of Theorem~\ref{thm:main}.
This is not implied by Kleiman's example~\cite{MR50:13063} of a subvariety in
a Grassmannian not meeting any translate of a particular Schubert variety
transversally, for his Schubert variety has codimension 2.
This emphasizes the need for new, perhaps less general or
more specific, criteria that imply transversality.

Let $X,Y,Z$, and $S$ be smooth equidimensional varieties 
with $\pi,f$, and $g$ the maps of Figure~\ref{fig:one},
\begin{figure}[htb]
$$
  \begin{picture}(94,99)(-2,2)
   \put(-5, 90){$X\times_Z Y$}  \put(80, 90){$Y$}
   \put( 8, 40){$X$}            \put(80, 40){$Z$}
   \put( 8,  0){$S$}
   \put(53, 98){$p_2$}
   \put(53,53){$f$}
   \put(-1,70){$p_1$}   \put(87,70){$g$}
   \put(1,23){$\pi$}
   \put( 13,85){\vector(0,-1){30}}
   \put( 83,85){\vector(0,-1){30}}
   \put( 13,35){\vector(0,-1){20}}
   \put( 38,93){\vector(1,0){35}}
   \put( 28, 43){\vector(1,0){45}}
  \end{picture}
$$
\caption{Intersection with a family\label{fig:one}}
\end{figure}
where $\pi$ is smooth, $g$ is unrammified, and $f$ is flat.
The fibres of $X\times_Z Y\to S$ are the intersections of the fibres of
$X\to S$ with $Y$ (along the maps $f$ and $g$).

Consider the bundles $E:=p_1^*T(X/S)$, the pullback of the relative tangent
bundle of $X\to S$ and $F:=p_2^*(g^*TZ)/TY$, the pullback of the conormal
bundle of $g$.
Let 
$$
  E \ \stackrel{\alpha}{\relbar\joinrel\longrightarrow}\ F
$$
be the map induced by $df\colon T(X/S)\to TZ$ and 
define $\Delta\subset X\times_Z Y$ to be the degeneracy locus of the
map $\alpha$, the set of points where the map $\alpha$ does not have full
rank.
This is in fact the locus where the intersection is not transverse.
Set
$$
  \begin{array}{rrl}
   \rho&:=& |\mbox{\rm rank}E - \mbox{\rm rank}F| +1\\
       &=& \dim X\times_Z Y -\dim S + 1\,. \qquad  \rule{0pt}{15pt}
           \mbox{\rm (when $X\times_Z Y$ dominates $S$)}
  \end{array}
$$
Then either $\Delta$ is empty or it has codimension 
at most $\rho$ in $X\times_Z Y$.

\begin{prop}[\cite{LaSp92}]
 Suppose either that $\Delta=\emptyset$ or else the
 codimension of $\Delta$ in $X\times_Z Y$ is equal to $\rho$.
 Then there is a non empty open subset $U$ of $S$ such that for $s\in U$
 either $X_s\times_ZY$ is empty, or it is smooth of the expected dimension 
 $\rho-1$.
\end{prop}

Laksov and Speiser define the family $X$ to be {\it determinantal}
if
\begin{equation}\label{determinantal}
 \begin{minipage}[c]{5.1in}
  \begin{center}
   \rule{0pt}{14pt}%
   For every unramified map $f\colon Y\to Z$ from a smooth variety $Y$,
   either\vspace{3pt}\\ $\Delta$ is empty, or else it has codimension 
   $\rho$ in $X\times_Z Y$, where $\rho$ is as above.\vspace{2pt}
  \end{center}
 \end{minipage}
\end{equation}
This condition is quite strong, 
as it implies that the general member of the family $X\to S$ meets
every unramified map $Y\to Z$ from a smooth variety $Y$
transversally.
Since enumerative geometry is concerned only with those maps
$Y\to Z$ that are `geometrically meaningful'---an admittedly ill-defined
class---this condition is perhaps stronger than needed by enumerative geometry.

They introduce a local condition which implies~\eqref{determinantal}.
Fix $z\in Z$ and a linear subspace $L\subset T_zZ$.
For any $x\in f^{-1}(z)$, we have the map
 \begin{equation}\label{alphaLx}
  T_x(X/S)\ \stackrel{\alpha_{L,x}}%
     {\relbar\joinrel\relbar\joinrel\relbar\joinrel\longrightarrow}\ T_zZ/L\,.
 \end{equation}
Set $\rho_L:=\dim X-\dim S +\dim L -\dim Z +1$
and let $\Delta_L\subset f^{-1}(z)$ be the locus of points $x$ where 
$\alpha_{L,x}$ has less than maximal rank.
Then either $\Delta_L=\emptyset$ or else it has codimension at most
$\rho_L$.
The family $X$ is {\it determinantal at $z$} if, for every
linear subspace $L\subset T_zZ$, either $\Delta_L=\emptyset$ or else it has
codimension $\rho_L$ in $f^{-1}(z)$.

\begin{prop}[\cite{LaSp92}]
If the family $X$ is determinantal at every $z\in Z$, then is it
determinantal.
\end{prop}\smallskip

Consider the enumerative problem (i) of \S 1 involving simple
Schubert varieties of ${\bf G}(r,n)$.
Let $l\geq 1$ be an integer, $Y:={\bf G}(r,n)$, and 
$Z:=[{\bf G}(r,n)]^l$, with the map $g\colon Y\to Z$ the diagonal embedding.
Let $S:=[{\bf G}(n{-}r,n)]^l$, and set
$$
  X\ :=\ \{(K_1,K_2,\ldots,K_l,\, H_1,H_2,\ldots,H_l)\in S\times Z\mid
            \dim K_i\cap H_i=1\}\,,
$$
which is the smooth points in the $l$-fold product of the universal
simple Schubert variety in ${\bf G}(n{-}r,n)\times{\bf G}(r,n)$
$$
  \{(K,H)\in{\bf G}(n{-}r,n)\times{\bf G}(r,n)\mid
       \dim K\cap H\neq\{0\}\}\,.
$$
In this case, the maps $\pi$ and $f$ are just the projections, and the fibre 
of $X\times_Y Z$ over a point $(K_1,K_2,\ldots,K_l)$ of $S$ is the intersection
 \begin{equation}\label{sm-fibre}
   \Omega^{\rm sm}(K_1)\cap\Omega^{\rm sm}(K_2)\cap
      \cdots\cap \Omega^{\rm sm}(K_l)\,,
\end{equation}
where $\Omega^{\rm sm}(K)$ consists of the smooth points of $\Omega(K)$.

Theorem~\ref{thm:main} asserts that when $l=r(n{-}r)$ (and hence for all 
$l\leq r(n{-}r)$), there is an open subset $U$ of $S$ consisting of points
$(K_1,K_2,\ldots,K_l)$ such that the intersection~(\ref{sm-fibre}) is non-empty
and is transverse at the generic point of each component.
This fact is not implied by the theory of Laksov and Speiser.

\begin{thm}\label{not-det}
   When $r, n{-}r>1$, the family $X$ is not determinantal at any $z\in Z$.
\end{thm}

Similar arguments show that transversality in the other enumerative problems
of Theorem~\ref{thm:main} cannot be obtained from the theory of Laksov and
Speiser.
For those, we replace $Y$ by one of 
${\mathbb F}\ell_{\bf r},OG(r)$, or ${\mathcal M}^q_{r,n}$,
and possibly modify $X,S$, and $Z$.\smallskip

Let $V=\K^n$.
For $H\in{\bf G}(r,n)$, the tangent space $T_H{\bf G}(r,n)$ equals
$\mbox{\rm Hom}(H,V/H)$.
Similarly, if $K\in{\bf G}(n{-}r,n)$ and $\dim H\cap K=1$, then 
(see, for instance~\cite[\S 2.9]{So97d})
$$
  T_H\Omega(K)\ =\ \{\varphi\in\mbox{\rm Hom}(H,V/H)\mid
          \varphi(K\cap H)\subset (H+K)/H\}\,.
$$
If we let $v=K\cap H\in{\mathbb P}(H)$ and 
$\Lambda=K{+}H\in{\mathbb P}^\vee(V/H)$, a hyperplane containing $H$,
then $T_H\Omega(K)=\tau_{v,\Lambda}$, where
$$
  \tau_{v,\Lambda}\ :=\ \{\varphi\in\mbox{\rm Hom}(H,V/H)\mid
          \varphi(v)\subset \Lambda/H\}\,.
$$
The key point 
is that the set of hyperplanes $T_H\Omega(K)$, for
$K\in\Omega^{\rm sm}(H)$, forms the proper subvariety of the space
${\mathbb P}^\vee(\mbox{\rm Hom}(H,V/H))$ consisting of hyperplanes
$\tau_{v,\Lambda}$ for 
$(v,\Lambda)\in{\mathbb P}(H)\times{\mathbb P}^\vee(V/H)$.
This also implies that the grassmannian 
hypothesis~\cite[Theorem 3.3]{LaSp90} fails to hold.
This set of hyperplanes is a single orbit of the stabilizer of $H$ in
$GL(V)$ of dimension $n{-}2$.
(Compare with the hypotheses of~\cite[Theorem 10]{MR50:13063}.)
\smallskip

\noindent{\sl Proof of Theorem~\ref{not-det}.}
Let $z=(K_1,K_2,\ldots,K_l)\in Z$.
Then
$$
  f^{-1}(z)\ =\  \{(K_1,K_2,\ldots,K_l)\in S\mid \dim K_i\cap H_i=1\}
      \  = \ \prod_{i=1}^l \Omega^{\rm sm}(H_i)\,,
$$
which has dimension $l[r(n{-}r){-}1]$.
Fix any $(v,\Lambda)\in{\mathbb P}(H)\times{\mathbb P}^\vee(V/H)$ and set
$$
  L\ :=\ \tau_{v,\Lambda}\times \prod_{i=2}^l T_{H_i}{\bf G}(r,n)\,.
$$

For $x=(K_1,K_2,\ldots,K_l)\in f^{-1}(z)$, we have
$$
  T_x(X/S)\ =\ \prod_{i=1}^l T_{H_i}\Omega(K_i)\ \subset\ 
               \prod_{i=1}^l T_{H_i}{\bf G}(r,n)\ =\ T_zZ\,,
$$
and the map $\alpha_{L,x}$~(\ref{alphaLx}) is the composition
$$
  T_x(X/S)\ \hookrightarrow\ T_zZ\ \twoheadrightarrow\ T_zZ/L\ 
      \simeq\ \K\,.
$$

This map drops rank precisely when $T_{H_1}\Omega(K_1)=\tau_{v,\Lambda}$,
that is, when $v\in K_1\subset \Lambda$ (and $v=K_1\cap H_1$).
This defines an open subset $\Omega'(v,\Lambda)$ of a Schubert subvariety of
${\bf G}(n{-}r,n)$ isomorphic to ${\bf G}(n{-}r{-}1,r{-}1)$, which has
dimension $(n{-}r{-}1)(r{-}1)=r(n{-}r){-}1{-}(n{-}2)$.
Thus
$$
  \Delta_L\ =\ \Omega'(v,\Lambda)\times
          \prod_{i=2}^l\Omega^{\rm sm}(K_i)
$$
has dimension $l[r(n{-}r){-}1]{-}(n{-}2)$, and hence
codimension $n{-}2$ in $f^{-1}(z)$.
However, 
 \begin{eqnarray*}
  \rho_L&=&\dim X-\dim S+\dim L-\dim Z + 1\\
        &=& l[r(n{-}r){-}1]\,,
 \end{eqnarray*}
which exceeds $n{-}2$ for all $l>0$ and $r,n{-}r>1$.
\qed


\end{document}